\documentclass[11pt, reqno]{amsart}

\usepackage[usenames,dvipsnames]{color}
\usepackage[OT2, T1]{fontenc}
\usepackage{url}
\usepackage{amsmath}
\usepackage{array}
\usepackage{graphicx}
\usepackage{amssymb}
\usepackage{amsthm}
\definecolor{link}{RGB}{11,0,128}
\usepackage[colorlinks=true,citecolor=NavyBlue,linkcolor=OliveGreen,urlcolor=link]{hyperref}
\usepackage{hyperref}
\usepackage{paralist}
\usepackage{colonequals}			% for nice :=
\usepackage{color}
\usepackage{mathabx}			% for \widec
\usepackage{amscd}
\usepackage[all,cmtip]{xy}			% for commutative diagrams
\usepackage{sseq}				% for spectral sequences
\usepackage{verbatim}
\usepackage{parskip}
\usepackage{microtype}
\usepackage[in]{fullpage}
\usepackage{mathrsfs} 			% for \mathscr (script letters)
\usepackage{cleveref}
\usepackage{enumitem}
\usepackage{moreenum}			% for enumeration via Greek letters
\usepackage[alphabetic,lite]{amsrefs} 	% for bibliography
\usepackage{stmaryrd}			% this is for \llbracket and \rrbracket
\usepackage[foot]{amsaddr}
\usepackage{subfiles}
\usepackage{ragged2e}
\usepackage{stackengine}

\newcommand{\bA}{\mathbb{A}}

\newcommand{\bC}{\mathbb{C}}

\newcommand{\bG}{\mathbb{G}}

\newcommand{\bP}{\mathbb{P}}

\newcommand{\bZ}{\mathbb{Z}}

\newcommand{\bbB}{\mathbf{B}}

\newcommand{\cE}{\mathcal{E}}

\newcommand{\cG}{\mathcal{G}}

\newcommand{\cP}{\mathcal{P}}

\newcommand{\cX}{\mathcal{X}}

\newcommand{\fp}{\mathfrak{p}}

\newcommand{\sO}{\mathscr{O}}

\DeclareMathOperator{\Aut}{Aut}		% The group of automorphisms
\DeclareMathOperator{\Bun}{Bun}			% The stack of bundles on a curve
\DeclareMathOperator{\GL}{GL}		% The general linear group
\DeclareMathOperator{\Gr}{Gr}			% Affine Grassmannian
\DeclareMathOperator{\Ker}{Ker}		% The kernel
\DeclareMathOperator{\Res}{Res}		% The restriction of the representation
\DeclareSymbolFont{cyrletters}{OT2}{wncyr}{m}{n}
\DeclareMathSymbol{\Sha}{\mathalpha}{cyrletters}{"58}	% Sha
\DeclareMathOperator{\Spec}{Spec}		% Spectrum of a ring
\newcommand{\ad}{\mathrm{ad}}			% associated adic space
\newcommand{\alg}{\mathrm{alg}}		% algebraic closure (mainly used in superscripts)

\newcommand{\ce}{\colonequals}
\newcommand{\dR}{{\mathrm{dR}}}		% de Rham

\newcommand{\et}{\mathrm{\acute{e}t}}	% for etale cohomology (mainly used in subscripts)
\newcommand{\fpqc}{\mathrm{fpqc}}		% for fpqc cohomology (mainly used in subscripts)
\newcommand{\gp}{{\mathrm{gp}}}		% Group completion of a monoid
\newcommand{\hra}{\hookrightarrow}
\renewcommand{\i}{^{-1}}
\newcommand{\isomto}{\overset{\sim}{\longrightarrow}}

\newcommand{\llb}{\llbracket}			% [[
\newcommand{\llp}{(\!(}			% ((
\newcommand{\ov}{\overline}
\providecommand{\p}[1]{\left(#1\right)}
\newcommand{\ra}{\rightarrow}

\newcommand{\rrb}{\rrbracket}			% ]]
\newcommand{\rrp}{)\!)}			% ))
\newcommand{\tensor}{\otimes} 			% tensor product
\newcommand{\Zar}{\mathrm{Zar}}		% for Zariski cohomology 

\providecommand{\up}[1]{{\upshape(}#1{\upshape)}}
\providecommand{\uref}[1]{{\upshape\ref{#1}}}
\providecommand{\uS}{{\upshape\S}}
\providecommand{\f}[2]{\frac{#1}{#2}}
\renewcommand{\b}{\textbf}
\providecommand{\ucolon}{{\upshape:} }

\newcommand{\brems}{\begin{rems} \hfill \begin{enumerate}[label=\b{\thenumberingbase.},ref=\thenumberingbase]}

\newcommand{\erems}{\end{enumerate} \end{rems}}
\newcommand{\begs}{\begin{egs} \hfill \begin{enumerate}[label=\b{\thenumberingbase.},ref=\thenumberingbase]}

\newcommand{\eegs}{\end{enumerate} \end{egs}}
\newcommand{\m}{\item}
\newcommand{\bsm}{\begin{smallmatrix}}
\newcommand{\esm}{\end{smallmatrix}}
\newcommand{\blem}{\begin{lemma}}
\newcommand{\elem}{\end{lemma}}

\newcommand{\bconj}{\begin{conj}}
\newcommand{\econj}{\end{conj}}
\newcommand{\bprob}{\begin{Problem}}
\newcommand{\eprob}{\end{Problem}}

\newcommand{\bq}{\begin{Q}}
\newcommand{\eq}{\end{Q}}
\newcommand{\benum}{\begin{enumerate}[label={{\upshape(\alph*)}}]}
\newcommand{\benuma}{\begin{enumerate}[label={{\upshape(\arabic*)}}]}
\newcommand{\benumb}{\begin{enumerate}[label={{\upshape\b{\arabic*.}}}]}
\newcommand{\benumr}{\begin{enumerate}[label={{\upshape(\roman*)}}]}
\newcommand{\eenum}{\end{enumerate}}
\newcommand{\bitem}{\begin{itemize}}
\newcommand{\eitem}{\end{itemize}}
\newcommand{\bc}{\begin{comment}}
\newcommand{\ec}{\end{comment}}
\newcommand{\bd}{\begin{defn}}
\newcommand{\ed}{\end{defn}}

\newcommand{\beg}{\begin{eg}}
\newcommand{\eeg}{\end{eg}}

\newcommand{\bcl}{\begin{claim}}
\newcommand{\ecl}{\end{claim}}

\newcommand{\x}{\text}

\newcommand{\q}{\quad}
\providecommand{\qxq}[1]{\quad\text{#1}\quad}
\providecommand{\qx}[1]{\quad\text{#1}}

\newcommand{\qq}{\quad\quad}
\newcommand{\qqq}{\quad\quad\quad}
\newcommand{\qqqq}{\quad\quad\quad\quad}

\newcommand{\tst}{\textstyle}
\newcommand{\ba}{\begin{aligned}}
\newcommand{\ea}{\end{aligned}}
\newcommand{\be}{\begin{equation}}
\newcommand{\ee}{\end{equation}}
\newcommand{\bpf}{\begin{proof}}
\newcommand{\epf}{\end{proof}}
\newcommand{\bthm}{\begin{thm}}
\newcommand{\ethm}{\end{thm}}
\newcommand{\bprop}{\begin{prop}}
\newcommand{\eprop}{\end{prop}}
\newcommand{\bcor}{\begin{cor}}
\newcommand{\ecor}{\end{cor}}
\newcommand{\brem}{\begin{rem}}
\newcommand{\erem}{\end{rem}}
\newcommand*{\QED}{\hfill\ensuremath{\qed}}

\usepackage{aliascnt}
\newaliascnt{numberingbase}{section}
\numberwithin{equation}{numberingbase}

\newtheoremstyle{thms}{0.5em}{0.5em}{\itshape}{}{\bfseries}{.}{ }{}
\theoremstyle{thms}
\newtheorem{conj}[numberingbase]{Conjecture}

\newtheorem{cor}[numberingbase]{Corollary}

\newtheorem{lemma}[numberingbase]{Lemma}

\newtheorem{prop}[numberingbase]{Proposition}

\newtheorem{Q}[numberingbase]{Question}
\newtheorem{thm}[numberingbase]{Theorem}

\newtheoremstyle{claims}{0.5em}{0.5em}{}{}{\itshape}{.}{ }{}
\theoremstyle{claims}
\newtheorem{claim}[equation]{Claim}

\newtheoremstyle{defs}{0.5em}{0.5em}{}{}{\bfseries}{.}{ }{}
\theoremstyle{defs}
\newtheorem{defn}[numberingbase]{Definition}

\newtheorem{eg}[numberingbase]{Example}
\newtheorem*{egs}{Examples}
\newtheorem{rem}[numberingbase]{Remark}

\newtheorem*{rems}{Remarks}

\Crefname{claim}{Claim}{Claims}
\Crefname{bclaim}{Claim}{Claims}
\Crefname{sublemma}{Lemma}{Lemmas}
\Crefname{conj}{Conjecture}{Conjectures}
\Crefname{cor}{Corollary}{Corollaries}
\Crefname{defn}{Definition}{Definitions}
\Crefname{eg}{Example}{Examples}
\Crefname{prop}{Proposition}{Propositions} 
\Crefname{Q}{Question}{Questions}
\Crefname{rem}{Remark}{Remarks}
\Crefname{thm}{Theorem}{Theorems}
\Crefname{Theorem}{Theorem}{Theorems}
\Crefname{variant}{Variant}{Variants}
\Crefname{caution}{Caution}{Cautions}

\theoremstyle{thms}
\newtheorem{thm-tweak}[subsection]{Theorem}
\Crefname{thm-tweak}{Theorem}{Theorems}
\newtheorem{lemma-tweak}[subsection]{Lemma}
\Crefname{lemma-tweak}{Lemma}{Lemmas}
\newtheorem{cor-tweak}[subsection]{Corollary}
\Crefname{cor-tweak}{Corollary}{Corollaries}
\newtheorem{prop-tweak}[subsection]{Proposition}
\Crefname{prop-tweak}{Proposition}{Propositions} 
\newtheorem{conj-tweak}[subsection]{Conjecture}
\Crefname{conj-tweak}{Conjecture}{Conjectures} 
\newtheorem{q-tweak}[subsection]{Question}
\Crefname{q-tweak}{Question}{Questions} 

\theoremstyle{defs}
\newtheorem{defn-tweak}[subsection]{Definition}
\Crefname{defn-tweak}{Definition}{Definitions}
\newtheorem{eg-tweak}[subsection]{Example}
\Crefname{eg-tweak}{Example}{Examples}
\newtheorem*{rems-tweak}{Remarks}
\newtheorem{rem-tweak}[subsection]{Remark}
\Crefname{rem-tweak}{Remark}{Remarks}

\newtheoremstyle{subsection-tweak}
   {2pt}
   {3pt}%
   {}
   {}%
   {\bfseries}
   {}%
   {.5em}
   {\thmnumber{\@{#1}{}\@{#2}.}%
    \thmnote{~{\bfseries#3.}}}    
    
\theoremstyle{subsection-tweak}
\newtheorem{pp}[numberingbase]{}
\newcommand{\bpp}{\begin{pp}}
\newcommand{\epp}{\end{pp}}

\theoremstyle{subsection-tweak}
\newtheorem{pp-tweak}[subsection]{}

\makeatletter
\def\@tocline#1#2#3#4#5#6#7{
    \begingroup 
    \@ifempty{#4}{}{}

    \parindent\z@ \leftskip#3\relax \advance\leftskip\@tempdima\relax
    #5\hskip-\@tempdima
      \ifcase #1
       \or\or \hskip 2em \or \hskip 1em \else \hskip 3em \fi%
      #6\nobreak\relax
    \dotfill\hbox to\@pnumwidth{\@tocpagenum{#7}}\par
    \nobreak
    \endgroup
 }
 \def\l@section{\@tocline{1}{0pt}{1pc}{}{}}

\renewcommand{\tocsection}[3]{%
  \indentlabel{\@ifnotempty{#2}{\makebox[1.3em][l]{%
    \ignorespaces#1 \bfseries{#2}.\hfill}}}\bfseries{#3}
    \vspace{-5pt}}

\renewcommand{\tocsubsection}[3]{%
  \indentlabel{\@ifnotempty{#2}{\hspace*{-0.5em}\makebox[2.1em][l]{%
    \ignorespaces#1#2.\hfill}}}#3
    \vspace{-5pt}}
   
\makeatother 

\makeatletter % This is for proper formatting of appendices in the table of contents
\newcommand\appendix@section[1]{%
  \refstepcounter{section}%
  \orig@section*{Appendix \@Alph\c@section. #1}%
}
\let\orig@section\section
\g@addto@macro\appendix{\let\section\appendix@section}
\makeatother

\makeatletter
\@namedef{subjclassname@2020}{%
  \textup{2020} Mathematics Subject Classification}
\makeatother

\author{K\k{e}stutis \v{C}esnavi\v{c}ius}
\address{CNRS, Universit\'{e} Paris-Saclay,   Laboratoire de math\'{e}matiques d'Orsay, F-91405, Orsay, France}
\email{kestutis@math.u-psud.fr}

\date{\today}

\begin{document}

\subjclass[2020]{Primary 14L15; Secondary 14M17.}
\keywords{Affine Grassmannian, loop group, reductive group, torsor.}

\title{The affine Grassmannian as a presheaf quotient}

\maketitle

\begin{abstract} 
For a reductive group $G$ over a ring $A$, its affine Grassmannian $\Gr_G$ plays important roles in a wide range of subjects and is typically defined as the \'etale sheafification of the presheaf quotient $LG/L^+G$ of the loop group $LG$ by its positive loop subgroup $L^+G$. We show that the Zariski sheafification gives the same result. Moreover, for totally isotropic $G$ (for instance, for quasi-split $G$), we show that no sheafification is needed at all: $\Gr_G$ is already the presheaf quotient $LG/L^+G$, which seems new already in the classical case of $G$ over $\bC$. 
For totally isotropic $G$, we also show that the affine Grassmannian may be formed using polynomial loops. We deduce all of these results from the study of $G$-torsors on $\bP^1_A$ that is ultimately built on the geometry of $\Bun_G$. 
\end{abstract}

\hypersetup{
    linktoc=page,     %set to all if you want both sections and subsections linked
}
\renewcommand*\contentsname{}
\q\\
\tableofcontents

The affine Grassmannian $\Gr_G$ of a reductive group $G$ originated in Lusztig's \cite{Lus83}*{Section 11} (see also Beilinson--Drinfeld's \cite{BD19}*{Section 4.5}) 
and is instrumental in the geometric Langlands program and other subjects that study $G$-torsors and their moduli. The goal of this article is to show that $\Gr_G$ admits a simpler definition than previously thought: for most $G$, it is simply the presheaf quotient $LG/L^+G$ of the loop group $LG$ by its positive loop subgroup $L^+G$, so that the  fpqc or \'etale sheafifications of this quotient that were used previously are not needed, see Theorems~\ref{thm:Zariski} and~\ref{thm:main} for precise statements.

\bpp[Conventions] \label{pp:conv}
As in \cite{SGA3IIInew}*{expos\'e XIX, d\'efinition 2.7}, a \emph{reductive group scheme} over a scheme $S$ is a smooth, affine $S$-group scheme whose geometric $S$-fibers are connected reductive groups. An $S$-scheme is \emph{locally constant} if fpqc locally on $S$ it becomes isomorphic to some $\bigsqcup_{i \in I} S$. 
\epp

\addtocounter{section}{-1}

\section{The affine Grassmannian and its modular description}

%\addtocounter{section}{1}

In this section, we fix a base ring $A$ and aim to review the definition and the modular interpretation of the affine Grassmannian $\Gr_\cG$ for a smooth, quasi-affine\footnote{Of course, over a field every quasi-affine group scheme of finite type is affine, see \cite{SGA3Inew}*{expos\'e VI$_{\x{B}}$, proposition~11.11}, and likewise for flat, finite type, separated groups over Dedekind rings with affine (or merely quasi-affine) generic fibers, see \cite{SGA3II}*{expos\'e XVII, proposition C.2.1 (3)} and \cite{Ana73}*{proposition 2.3.1}. Over higher-dimensional base rings, however, there exist quasi-affine groups that are not affine, see \cite{Ray70b}*{chapitre~VII, section 3} for such an example over $\bC[x, y]$.} $A\llb t\rrb$-group scheme $\cG$, see \Cref{prop:modular}.

\bpp[The loop functor] \label{pp:loops}
For a functor $\cX$ on the category of $A\llp t \rrp$-algebras (resp.,~$A\llb t \rrb$-algebras), its \emph{loop functor} $L\cX$ (resp.,~\emph{positive loop functor} $L^+\cX$) is defined on the category of $A$-algebras~by 
\[
L\cX\colon B \mapsto \cX(B\llp t \rrp) \qxq{(resp.,} L^+\cX\colon B \mapsto \cX(B\llb t \rrb)).
\]
The morphism $L^+\cX \ra L\cX$ is often an inclusion, for instance, this is so whenever $\cX$ is a subfunctor of a separated $A\llb t \rrb$-scheme (see \cite{EGAI}*{corollaire 9.5.6}). 

If $\cX$ is an $A\llb t\rrb$-scheme, then, by \cite{Bha16}*{Theorem 4.1 and Remark 4.6}, $L^+\cX$ is an $A$-scheme because
\[
\tst L^+\cX \cong \varprojlim_{n > 0} \Res_{(A[t]/(t^n))/A} (\cX_{A[t]/(t^n)}), \qxq{equivalently,} (L^+\cX)(B) = \varprojlim_{n > 0} \cX(B[t]/(t^n))
\]
for every $A$-algebra $B$, compare with \cite{CLNS18}*{Chapter 3, Corollary 3.3.7 b)}. Similarly, if $\cX$ is a quasi-compact and quasi-separated $A\llb t\rrb$-algebraic space, then, by \cite{Bha16}*{Theorem 4.1} and \cite{SP}*{Proposition~\href{https://stacks.math.columbia.edu/tag/05YF}{05YF} and Lemmas~\href{https://stacks.math.columbia.edu/tag/07SF}{07SF} and \href{https://stacks.math.columbia.edu/tag/05YD}{05YD}}, $L^+\cX$ is a quasi-compact and quasi-separated $A$-algebraic space. If $\cX$ is even an affine $A\llb t \rrb$-scheme, then, by considering a presentation of its coordinate ring in terms of generators and relations, $L^+\cX$ is an affine $A$-scheme and $L\cX$ is an ind-affine $A$-ind-scheme, more precisely, there are affine $A$-schemes $X_n$ and closed immersions
\[
\tst X_0 \hra X_1 \hra \dotsb \qxq{such that} L^+\cX = X_0 \qxq{and} L\cX = \bigcup_{n \ge 0} X_n \qx{as functors on $A$-algebras $B$.}
\]
For an $A$-algebra $B$, let $B\{t\}$ be the Henselization of $B[t]$ with respect to the ideal $tB[t]$, see \cite{Hitchin-torsors}*{Section 2.1.2} or  \cite{SP}*{Lemma \href{https://stacks.math.columbia.edu/tag/0A02}{0A02}}. It is useful to consider Henselian (resp.,~algebraic; resp.,~polynomial) variants $L_h \cX$ and $L^+_h\cX$ (resp.,~$L_\alg \cX$ and $L^+_\alg\cX$; resp.,~$L_{\mathrm{poly}} \cX$ and $L^+_{\mathrm{poly}}\cX$):
\[ \ba
 &\tst L_h \cX\colon B \mapsto \cX(B\{ t \}[\f1t]) \qqqq\ \ \,  \qxq{and} \q \tst  L^+_h\cX \colon B \mapsto \cX(B\{ t \}), \\
\tst & L_\alg\cX \colon B \mapsto \tst\cX((B[t]_{1 + tB[t]})[\f1t]) \q \qxq{and} \q \tst L^+_\alg\cX \colon B \mapsto \cX(B[t]_{1 + tB[t]}), \\
&\tst L_{\mathrm{poly}}\cX \colon B \mapsto \tst\cX(B[t, t\i]) \qqq \ \,  \qxq{and} \q L^+_{\mathrm{poly}}\cX \colon B \mapsto \cX(B[t]),
\ea \]
granted that $\cX$ is begins its life over $A\{t\}[\f1t]$ (resp.,~over $(A[t]_{1 + tA[t]})[\f1t]$; resp.,~over $A[t, t\i]$) and, for $L^+_*\cX$, even already over $A\{t\}$ (resp.,~over $A[t]_{1 + tA[t]}$; resp.,~over $A[t]$). These variant functors are sometimes easier to handle, for instance, they all commute with filtered direct limits in $B$ granted that so does $\cX$, %(in other words, granted that $\cX$ itself is locally of finite presentation), 
moreover, $L^+_{\mathrm{poly}}\cX$ and $L_{\mathrm{poly}}\cX$ are nothing else but restrictions of scalars.
%All of these functors $L^+_*\cG$ and $L_*\cG$ above are sheaves for the fpqc topology on $A$-algebras $B$. Indeed, for an  

By the following proposition, for many $\cX$, the functors $L_*^{(+)}\cX$ cannot be improved by sheafifying. 
\epp

\bprop \label{prop:L-sheaf}
For a scheme $\cX$ over $A\llp t \rrp$, or over $A\{t\}[\f1t]$, or over $(A[t]_{1 + tA[t]})[\f1t]$, or over $A[t, t\i]$ \up{resp.,~over $A\llb t \rrb$, or over $A\{t\}$, or over $A[t]_{1 + tA[t]}$, or over $A[t]$} as in \uS\uref{pp:loops} such that every quasi-compact open of $\cX$ is quasi-affine, the functor $L_*\cX$ \up{resp.,~and also its subfunctor $L^+_*\cX$} is an fpqc sheaf on the category of $A$-algebras.
\eprop

\bpf
Since $\cX$ is separated, \cite{EGAI}*{corollaire 9.5.6} ensures that the map $L^+_*\cX \ra L_*\cX$ is indeed an inclusion. Thus, all we need to show is that for an fpqc cover $B \ra B'$ of $A$-algebras, the sequence
\[
(L_*\cX)(B) \ra (L_*\cX)(B') \rightrightarrows (L_*\cX)(B' \tensor_B B') \ \ \x{(resp.,}\ \  (L_*^+\cX)(B) \ra (L_*^+\cX)(B') \rightrightarrows (L_*^+\cX)(B' \tensor_B B'))
\]
is exact. The case of the polynomial variant follows from fpqc descent \cite{SP}*{Lemma~\href{https://stacks.math.columbia.edu/tag/023Q}{023Q}} because
\[
B[t] \ra B'[t] \rightrightarrows  (B' \tensor_B B')[t] \qxq{and} B[t, t\i] \ra B'[t, t\i] \rightrightarrows  (B' \tensor_B B')[t, t\i]
\]
are both fpqc cover sequences. For other variants, since each ring-valued point of $\cX$ factors through a quasi-compact open, we may assume that $\cX$ is quasi-affine. The maps  
\[
B\{t\} \ra B'\{t\} \qxq{and} B[t]_{1 + tB[t]} \ra B'[t]_{1 + tB'[t]}
\]
 are both faithfully flat, in particular, they induce surjections on spectra, so the Henselian and algebraic variants reduce further to affine $\cX$. In the power series case, the reduction to affine $\cX$ is more subtle if $B$ and $B'$ are not Noetherian (because then we do not know whether the map $\Spec B'\llp t \rrp \ra \Spec B \llp t \rrp$ is surjective) and will use ideas from \cite{Hitchin-torsors}*{Lemma~2.2.9~(i)} as follows.

By realizing $\cX$ as the complement of the vanishing locus of a finitely generated ideal in an affine $A\llp t \rrp$-scheme (resp.,~$A\llb t\rrb$-scheme), to reduce to affine $\cX$ we need to show that elements $b_1, \dotsc, b_n$ in $B\llp t \rrp$ (resp.,~in $B\llb t \rrb$) generate the unit ideal as soon as their images do so in $B'\llp t \rrp$ (resp.,~in  $B'\llb t \rrb$). The case of $B\llb t \rrb$ follows by checking modulo $t$ and using the faithful flatness of $B \ra B'$. To treat $B\llp t \rrp$, we may assume that $b_1, \dotsc, b_n \in B\llb t \rrb$ and use the faithfully flat cover $B\llb t \rrb \ra (B\llb t \rrb \tensor_B B')^h_{(t)}$ (Henselization with respect to the ideal generated by $t$) to reduce to showing that $b_1, \dotsc, b_n$ generate the unit ideal in $(B\llb t \rrb \tensor_B B')^h_{(t)}[\f1t]$. The $t$-adic completion of $(B\llb t \rrb \tensor_B B')^h_{(t)}$ is $B'\llb t \rrb$ and, by assumption, there are $a_1, \dotsc, a_n \in B'\llb t \rrb$ with $a_1b_1 + \dotsc + a_nb_n = t^N$ in $B'\llb t \rrb$ for some $N > 0$. By approximating modulo $t^{N + 1}$, therefore, there are $a'_1, \dotsc, a'_n \in (B\llb t \rrb \tensor_B B')^h_{(t)}$ such that 
\[
a_1'b_1 + \dotsc + a_n'b_n \in t^N + t^{N + 1}(B\llb t \rrb \tensor_B B')^h_{(t)}.
\]
This means that $b_1, \dotsc, b_n$ indeed generate the unit ideal in $(B\llb t \rrb \tensor_B B')^h_{(t)}[\f1t]$, as desired. 

Now that $\cX$ is affine, the functor $\cX(-)$ turns fiber products of rings into fiber products of sets. Thus, all we need to show is the exactness of the horizontal sequences in the commutative diagram
\[\xymatrix@R=12pt{
B[t]_{1 + tB[t]} \ar@{^(->}[d] \ar[r] & B'[t]_{1 + tB'[t]}  \ar@{^(->}[d]\ar@<0.5ex>[r]\ar@<-0.5ex>[r] & (B' \tensor_B B')[t]_{1 + t(B' \tensor_B B')[t]} \ar@{^(->}[d] \\
B\{t\} \ar[r]\ar@{^(->}[d] &B'\{t\} \ar@<0.5ex>[r]\ar@<-0.5ex>[r]\ar@{^(->}[d] &(B' \tensor_B B')\{t\}\ar@{^(->}[d] \\
B\llb t\rrb \ar[r] &B'\llb t\rrb \ar@<0.5ex>[r]\ar@<-0.5ex>[r]  &(B' \tensor_B B')\llb t\rrb,
}
\]
which implies the corresponding exactness after further inverting $t$; here the bottom vertical maps are injective by \cite{Hitchin-torsors}*{Section 2.1.2} (a limit argument to reduce to finite type $\bZ$-algebras). The exactness of the bottom row is seen coefficientwise. Thus, by fpqc descent, it is enough to show that
\[\ba
B'\{t\} \tensor_{B\{t\}} B'\{t\} &\hra (B' \tensor_B B')\{ t\}, \\
B'[t]_{1 + tB'[t]} \tensor_{B[t]_{1 + tB[t]}} B'[t]_{1 + tB'[t]} &\hra (B' \tensor_B B')[t]_{1 + t(B' \tensor_B B')[t]}.
\ea\]
The injectivity of the second map is evident because $B'[t] \tensor_{B[t]} B'[t] \isomto (B'\tensor_B B^\prime)[t]$ and the elements of $1 + t(B'\tensor_B B')[t]$ are nonzerodivisors. As for the first map, if the $B$-algebra $B'$ was finitely presented, then we could use a limit argument to reduce to a Noetherian situation and then check the injectivity after passing to completions. In general, since every \'{e}tale $(B'\tensor_B B')[t]$-algebra $(B'\tensor_B B')$-fiberwise has no embedded associated primes, both of the maps in question are injections by \cite{RG71}*{premi\`ere partie, corollaire 3.2.6} and a limit argument.
\epf

\brem \label{rem:ind-quasi-affine}
In \Cref{prop:L-sheaf}, the assumption on quasi-compact opens  holds if $\cX$ is a subscheme of an affine scheme. In general, a scheme whose quasi-compact opens are all quasi-affine is \emph{ind-quasi-affine}, see \cite{Hitchin-torsors}*{Definition~2.2.5} and the references to \cite{SP} given there for details. We avoid this terminology here because it may be confusing in the context of ind-schemes. The automorphism scheme of a reductive group is such an $\cX$ that need not be quasi-affine, see \cite{torsors-regular}*{Section~1.3.7}. Another useful example is locally constant schemes $\cX$, for which even $L\cX = L^+\cX$ as follows. 

\erem

\bcor \label{cor:L-Lplus}
Let $A$ be a ring.
\benum
\m \label{m:LL-a}
In the following diagram, the squares are Cartesian and the maps are bijective on idempotents\ucolon 
\[
\q\xymatrix@C=15pt@R=15pt{
A \ar@{^(->}[r]\ar@{_(->}[rd] &A[t] \ar@{^(->}[r] \ar@{^(->}[d] & A[t]_{1 + tA[t]} \ar@{^(->}[r]\ar@{^(->}[d] & A\{t \} \ar@{^(->}[r]\ar@{^(->}[d] & A\llb t \rrb \ar@{^(->}[d] \\
& A[t, t\i] \ar@{^(->}[r]& (A[t]_{1 + tA[t]})[\f1t] \ar@{^(->}[r] &A\{t \}[\f1t] \ar@{^(->}[r] & A\llp t \rrp.
}
\]

\m \label{m:LL-b}
For every scheme $\cX$ over $A\llb t \rrb$, or over $A\{t\}$, or over $A[t]_{1 + tA[t]}$, or over $A[t]$ as in \uS\uref{pp:loops} such that $\cX$ is locally constant \up{see \uS\uref{pp:conv}}, we have $L_*^+\cX \isomto L_*\cX$.
\eenum
\ecor

\bpf\hfill
\benum
\m
As indicated, the maps are all injective, either by inspection, or by faithful flatness, or by a limit argument given in \cite{Hitchin-torsors}*{Section 2.1.2}. The squares are Cartesian by \cite{SP}*{Lemma~\href{https://stacks.math.columbia.edu/tag/0BNR}{0BNR}}. As for idempotents, since the maps are injective, it suffices to recall from \cite{Hitchin-torsors}*{Corollary~2.1.19} that the map $A \ra A\llp t \rrp$ is bijective on idempotents (as one may also show directly).

\m
By \cite{SP}*{Lemma \href{https://stacks.math.columbia.edu/tag/0AP8}{0AP8}}, we may check fpqc locally on the base that our locally constant $\cX$ satisfies the assumptions of \Cref{prop:L-sheaf}, more precisely, that every quasi-compact open of $\cX$ is quasi-affine. In particular, in the case when $\cX$ begins life over $A\{t\}$, the square
\[
\qq\xymatrix{
\cX(A\{t\}) \ar@{^(->}[r]\ar[d] & \cX(A\{t\}[\f1t]) \ar[d] \\
\cX(A\llb t \rrb) \ar@{^(->}[r] & \cX(A\llp t \rrp)
}
\]
is Cartesian by \cite{Hitchin-torsors}*{Proposition 2.2.12}, and likewise with $A[t]_{1 + tA[t]}$ or $A[t]$ in place of $A\{t\}$. Since we may vary $A$, this reduces us to showing that the inclusion $\cX(A\llb t\rrb) \subset \cX(A\llp t\rrp)$ is an equality. For this, since \emph{loc.~cit.}~also implies that for every prime $\fp \subset A$ the square  
\[
\qq\xymatrix{
\cX(A\llb t \rrb \tensor_A A_\fp) \ar@{^(->}[r]\ar[d] & \cX(A\llp t \rrp \tensor_A A_\fp) \ar[d] \\
\cX(A_\fp\llb t \rrb) \ar@{^(->}[r] & \cX(A_\fp\llp t \rrp)
}
\]
is Cartesian, by combining this square with spreading out, we may assume further that $A$ is local. At this point, \Cref{prop:L-sheaf} allows us replace $A$ by its strict Henselization.

Once $A$ is strictly Henselian local, we fix any $A\llp t \rrp$-point of $\cX$ and form its schematic image to get a closed subscheme $Z \subset \cX$, see \cite{SP}*{Definition  \href{https://stacks.math.columbia.edu/tag/01R7}{01R7}}. It is enough to show that 
\[
\qq Z \overset{?}{\cong} \Spec(A \llb t \rrb).
\]
We first claim that the map $Z \ra  \Spec(A \llb t \rrb)$ is surjective. For this, since $\cX$ is locally constant, we may choose a faithfully flat $A \llb t \rrb$-algebra $B$ such that $\cX_B \cong \bigsqcup_{i \in I}  \Spec B$. By \cite{SP}*{Lemma~\href{https://stacks.math.columbia.edu/tag/081I}{081I}}, the formation of $Z$ commutes with base change to $B$, so we need to show that the schematic image of any $B[\f1t]$-point of $\cX_B$ surjects onto $\Spec B$. It does because otherwise, by the explicit nature of $\cX_B$, the open immersion $\Spec (B[\f1t]) \subset \Spec B$ would factor through some proper closed subscheme, contradicting its schematic density. %inherited from $\Spec(A\llp t \rrp) \subset \Spec(A\llb t\rrb)$. 

Since $A\llb t \rrb$ is strictly Henselian local and $\cX$ is \'etale, \cite{EGAIV4}*{th\'eor\`eme~18.5.11} ensures that for every $z \in Z$ above the maximal ideal of $A\llb t \rrb$, we have 
\[
\qq \cX \cong \Spec (\sO_{\cX,\, z}) \sqcup \cX' \qxq{with} \sO_{\cX,\, z} \cong A\llb t \rrb.
\]
We have already argued that such a $z$ exists, and \ref{m:LL-a} shows that $\Spec(A\llp t \rrp)$ inherits connectedness from $\Spec A$, and so ensures that $Z \cap \cX' = \emptyset$. The desired $Z \cong \Spec(A \llb t \rrb)$ now follows from the schematic density of the open immersion $\Spec(A\llp t \rrp) \subset \Spec(A\llb t \rrb)$.  \qedhere
\eenum
\epf

\bpp[The affine Grassmannian] \label{pp:GrG}
In the setting of \S\ref{pp:loops}, if the functor $\cX$  is group valued, then so are the functors $L^+_*\cX$ and $L_*\cX$. With this in mind, for a smooth, quasi-affine $A\llb t \rrb$-group scheme $\cG$, its \emph{affine Grassmannian} is the pointed set valued functor defined by 
\be \label{eqn:GrG}
\Gr_\cG \ce (L\cG/L^+\cG)_\et,
\ee
where, as indicated, the sheafification of the presheaf quotient $L\cG/L^+\cG$ is formed in the \'etale topology. If instead $\cG$ begins life as a smooth, quasi-affine $A\{t\}$-group scheme, then \cite{Hitchin-torsors}*{Example 2.2.19} (which is based on approximation and algebraization techniques) ensures that we may use Henselian loops instead, more precisely, that we have an identification of presheaf quotients 
\[
L\cG/L^+\cG \cong L_h\cG/L_h^+\cG
\]
In many situations, the affine Grassmannian $\Gr_\cG$ is represented by an $A$-ind-scheme, which is often even ind-projective over $A$, see, for instance, \cite{PR08}. This proceeds by embedding $\cG$ into some $\GL_n$, which is not always possible with our general assumptions, so we will not use representability results. 
\epp

In the literature one sometimes finds the affine Grassmannian defined as the fpqc sheafification $(L\cG/L^+\cG)_{\mathrm{fpqc}}$. A priori this makes no mathematical sense: even for a field, isomorphism classes of its fpqc covers form a proper class, so the fpqc sheafification may not exist;  however, by the following proposition, this ``fpqc sheafification approach'' gives the same result for our $\cG$ as above.

\bprop \label{prop:modular}
For a ring $A$ and a smooth, quasi-affine $A\llb t \rrb$-group scheme $\cG$, the affine Grassmannian $\Gr_\cG$ is an fpqc sheaf and has the following modular description on the category of $A$-algebras~$B$\ucolon
\[
\Gr_\cG(B) \cong \left\{ (\cE, \iota) \,\, \, \Big\vert\,  \begin{aligned}&\  \cE\text{ is a }\cG\text{-torsor over }B\llb t\rrb,\\ & \  \iota \in  \cE(B\llp t \rrp) \text{ is a trivialization over $B\llp t \rrp$}
    \end{aligned}  \right\} /\sim,
\]
while the presheaf quotient $L\cG/L^+\cG$ is a subfunctor of $\Gr_\cG$ parametrizing those $(\cE, \iota)$ with $\cE$ trivial. 
\eprop

\bpf
For the trivial $\cG$-torsor, $L^+\cG$ parametrizes its $\cG$-torsor automorphisms over $(-)\llb t\rrb$ while $L\cG$ parametrizes its trivializations $\iota$ over $(-)\llp t \rrp$, so the presheaf quotient $L\cG/L^+\cG$ is the claimed subfunctor of $\Gr_\cG$, granted that the latter has the displayed modular description. Moreover, by the smoothness and quasi-affineness of $\cE$ inherited from $\cG$ (see \cite{SP}*{Lemma \href{https://stacks.math.columbia.edu/tag/0247}{0247}}) and by the infinitesimal lifting of sections (compare also with \cite{Hitchin-torsors}*{Proposition~2.1.4}), each $(\cE, \iota)$ lands in this subfunctor \'{e}tale locally on $B$, so the sought modular description will follow once we argue that it defines an fpqc sheaf. For this, since the pairs $(\cE, \iota)$ have no nontrivial automorphisms, all that remains is to show that the groupoid-valued~functor
\[
B \mapsto \left\{ (\cE, \iota) \,\, \, \Big\vert\,  \begin{aligned}&\  \cE\text{ is a }\cG\text{-torsor over }B\llb t\rrb,\\ & \  \iota \in  \cE(B\llp t \rrp) \text{ is a trivialization over $B\llp t \rrp$}
    \end{aligned}  \right\}
\]
is a stack for the fpqc topology on $A$-algebras $B$. For this, let $B \ra B'$ be an fpqc cover of $A$-algebras and let $(\cE', \iota')$ be a pair for $B'$ equipped with a descent datum with respect to this cover. It suffices to uniquely descend $\cE'$ to a $\cG$-torsor $\cE$ over $B\llb t \rrb$: the trivialization $\iota'$ will then also descend because $L\cE$ is an fpqc sheaf by \Cref{prop:L-sheaf}. The usual fpqc descent for $\cG$-torsors (see \cite{SP}*{Lemma \href{https://stacks.math.columbia.edu/tag/0247}{0247}}) gives us a unique compatible system of $\cG$-torsor descents $\cE_n$ over $B\llb t\rrb/(t^{n+1})$ for $n \ge 0$. It then remains to uniquely algebraize this sequence to a $\cG$-torsor $\cE$ over $B\llb t \rrb$, and there are several ways to do this (with even more possible arguments if $\cG$ is affine). Perhaps the most direct is to apply the algebraization result \cite{BHL17}*{Theorem 8.1} to the classifying stack $\bbB \cG$. A somewhat more elementary approach is to first use \cite{Hitchin-torsors}*{Theorem~2.1.6} to lift $\cE_0$ to a $\cG$-torsor $\cE$ over $B\llb t \rrb$ and to then inductively build a compatible sequence of isomorphisms 
\[
\cE|_{B\llb t\rrb/(t^{n + 1})} \simeq \cE_n,
\]
by using \emph{loc.~cit.}~and the smoothness of the quasi-affine $B\llb t\rrb$-group $\underline{\Aut}_\cG(\cE)$. 
\epf

\addtocounter{section}{-6}

\section{The affine Grassmannian as a  Zariski quotient}

\addtocounter{section}{5}

Even though the affine Grassmannian of $\cG$ is defined as the \'etale sheafification of the presheaf quotient $L\cG/L\cG^+$ (see \eqref{eqn:GrG}), in \Cref{thm:Zariski} below we show that the Zariski sheafification suffices when $\cG$ is reductive and descends to $A$. For intuition for why the Zariski topology may be enough, we recall that the inclusion $L^+\cG \subset L\cG$ is morally similar to the inclusion $P \subset G$ of a parabolic subgroup of a reductive group scheme, and that $(G/P)_\fpqc \cong (G/P)_\Zar$ because over any semilocal ring parabolic subgroups of the same type are conjugate, see \cite{SGA3IIInew}*{expos\'e XXVI, corollaire~5.2}. 

We will deduce that the Zariski sheafification is enough  from the following result about torsors over $\bP^1_A$ that is proved by studying the geometry of the algebraic stack $\Bun_G$ that parametrizes such torsors, or in \cite{PS23b} by a different approach. Various weaker and more technical earlier variants of this result would suffice as well, for instance, \cite{Fed22a}*{Theorem~6} or \cite{torsors-regular}*{Proposition~5.3.6}. In effect, in some sense, we obtain the main results of this article by ascending geometric information along the uniformization map $\Gr_G \ra \Bun_G$.

\bthm[\cite{totally-isotropic}*{Theorem 3.6}] \label{thm:CF-input}
Let $G$ be a reductive group scheme over a semilocal ring $A$. Every $G$-torsor $E$ over $\bP^1_A$ is $A$-sectionwise constant, equivalently, $E|_{\{ t= 0 \}} \simeq E|_{\{ t= \infty \}}$. \QED
\ethm

We will use \Cref{thm:CF-input} through its following consequence for reductive group torsors over $A\llb t \rrb$. 

\bprop \label{prop:semilocally}
Let $A$ be a semilocal ring and let $\cG$ be an $A\llb t \rrb$-group scheme that is an extension of an $A\llb t \rrb$-group $\cG^\et$ that is locally constant \up{see \uS\uref{pp:conv}} by a reductive $A\llb t \rrb$-group scheme $\cG^0$ such that $\cG^0_{A\llp t \rrp}$ descends to a reductive $A$-group scheme. No nontrivial $\cG$-torsor over $A\llb t \rrb$ trivializes over $A\llp t \rrp$, that is, we have
\[
\Ker\left(H^1(A\llb t \rrb, \cG) \ra H^1(A\llp t \rrp, \cG)\right) = \{ *\};
\]
in particular, $\cG^0$ itself descends to a reductive $A$-group scheme. 
\eprop

\bpf
Let $\cE$ be a $\cG$-torsor over $A\llb t \rrb$ that trivializes over $A\llp t \rrp$. By \Cref{cor:L-Lplus}~\ref{m:LL-b}, every $e \in \cE(A\llp t\rrp)$ gives rise to an $A\llp t \rrp$-point of the $\cG^\et$-torsor $\cE/\cG^0$ that extends uniquely to an $A\llb t \rrb$-point $\ov{e} \in (\cE/\cG^0)(A\llb t \rrb)$. The preimage of $\ov{e}$ in $\cE$ is a $\cG^0$-torsor that trivializes over $A\llp t \rrp$. We may replace $\cE$ by this preimage to reduce to the case when $\cG = \cG^0$. 

Now that $\cG = \cG^0$, suppose first that $\cG$ descends to a reductive $A$-group $G$. Due to the triviality over $A\llp t \rrp$, we may patch $\cE$ with the trivial $G$-torsor over $\bP^1_A \setminus \{ t = 0\}$ (see, for instance, \cite{Hitchin-torsors}*{Lemma~2.2.11~(b)}) to build a $G$-torsor $E$ over $\bP^1_A$ such that $E|_{\{t = 0\}} \simeq \cE|_{\{t = 0\}}$ and $E|_{\{t = \infty\}}$ is trivial. By \Cref{thm:CF-input}, then $\cE|_{\{t = 0\}}$ is also trivial, to the effect that, by the infinitesimal lifting of sections due to smoothness, $\cE$ is trivial, too, as desired (compare also with \cite{Hitchin-torsors}*{Proposition 2.1.4}).

To complete the proof, it remains to show that the reductive $A$-group scheme $G$ for which we have $G_{A\llp t \rrp} \simeq \cG_{A\llp t \rrp}$ also descends $\cG$, that is, that already $G_{A\llb t \rrb} \simeq \cG_{A\llb t \rrb}$. By \Cref{cor:L-Lplus}~\ref{m:LL-a} and the classification of reductive group schemes (recalled in \cite{torsors-regular}*{Sections 1.3.1 and 1.3.7}), our $\cG$ corresponds to an $\underline{\Aut}_\gp(G)$-torsor over $A\llb t \rrb$ that trivializes over $A\llp t \rrp$. However, as recalled in \emph{loc.~cit.},~$\underline{\Aut}_\gp(G)$ is an extension of a locally constant $A$-group by the reductive $A$-group $G^\ad$. In particular, the case settled in the first two paragraphs of the proof applies and shows that the $\underline{\Aut}_\gp(G)$-torsor in question is trivial already over $A\llb t \rrb$, so that $G_{A\llb t \rrb}\simeq \cG$, as desired. 
\epf

\bcor \label{cor:H1-local}
For a reductive group scheme $G$ over a semilocal ring $A$, we have %no nontrivial $G$-torsor trivializes over $A\llp t \rrp$, in other words, we have
\[
H^1(A, G) \hra H^1(A\llp t \rrp, G).
\]
\ecor

\bpf
By twisting (see \cite{torsors-regular}*{equation (1.2.1.1)}), it suffices to show that the map in question has trivial kernel. 
This follows from \Cref{prop:semilocally} because, by \cite{Hitchin-torsors}*{Theorem 2.1.6}, we have 
\[
H^1(A, G) \isomto H^1(A\llb t \rrb, G). \qedhere
\] 
\epf

\brem
Results of \cite{FG21} suggest that the reductivity assumption may be nonessential for \Cref{cor:H1-local}. It would therefore be interesting to find a more general result of this type. 
\erem

We are ready for the promised sufficiency of the Zariski sheafification for $\Gr_\cG$ in the case when $\cG$ is (constant) reductive. One may compare this to its earlier variant \cite{Bac19}*{Proposition 14} that restricted to smooth $A$ over a field and deduced the conclusion from the Grothendieck--Serre conjecture. For us, it is the study of the latter that has indirectly led to the results of this article.  

\bthm \label{thm:Zariski}
For a ring $A$ and an $A\llb t \rrb$-group scheme $\cG$ that is an extension of a finite \'etale $A\llb t \rrb$-group $\cG^\et$ by a reductive $A\llb t \rrb$-group scheme $\cG^0$ that descends to a reductive $A$-group scheme, the affine Grassmannian $\Gr_\cG$ is the Zariski sheafification of the presheaf quotient $L\cG/L^+\cG$, that is,
\[
\Gr_\cG \cong (L\cG/L^+\cG)_\Zar.
\]
\ethm

\bpf
By the modular interpretation supplied by \Cref{prop:modular}, the possibility of varying $A$, and the infinitesimal lifting of sections due to smoothness, all we need to show is that for any $\cG$-torsor $\cE$ over $A\llb t \rrb$ that trivializes over $A\llp t \rrp$, the $\cG$-torsor $\cE|_{\{t = 0\}}$ trivializes Zariski locally on $A$. However, \Cref{prop:semilocally} ensures that $\cE|_{\{t = 0\}}$ trivializes even Zariski semilocally on $A$. 
\epf

\brem \hfill
\benuma
\m
We do not know whether the reductivity assumption is critical for \Cref{thm:Zariski}, for instance, whether the Zariski sheafification also suffices for parahoric groups. Certainly, it does in the case when $\cG$ is a reductive $A\llb t \rrb$-group scheme and $\cP$ is a smooth, quasi-affine $A\llb t \rrb$-group scheme equipped with an $A\llb t \rrb$-morphism $\cP \ra \cG$ that modulo $t$ reduces to an inclusion of a parabolic subgroup: indeed, by \cite{Hitchin-torsors}*{Theorem 2.1.6} and \cite{torsors-regular}*{equation (1.3.5.2)}, these assumptions ensure that $H^1(B\llb t \rrb, \cP) \subset H^1(B \llb t \rrb, \cG)$ for any semilocal $A$-algebra $B$, so 
\[
\qqq\Ker\p{H^1(B\llb t \rrb, \cG) \ra H^1(B\llp t \rrp, \cG) } = \{ *\} \  \Longrightarrow \ \Ker\p{H^1(B\llb t \rrb, \cP) \ra H^1(B\llp t \rrp, \cP) } = \{ *\},
\]
to the effect that \Cref{prop:modular} and \Cref{thm:Zariski} imply the claimed $\Gr_\cP \cong (L\cP/L^+\cP)_\Zar$. 

\m
We do not know whether \Cref{thm:Zariski} admits a version for the Witt vector affine Grassmannian. For a version of \Cref{thm:Zariski} for the $B_\dR^+$-affine Grassmannian, see \cite{BdR-sheafification}*{Theorem~3.1}.
\eenum
\erem

%\bpp[Notation and conventions] \label{conv}

%\epp

\addtocounter{section}{-11}

\section{The affine Grassmannian as a presheaf quotient}

\addtocounter{section}{10}

For most reductive group schemes $\cG$, even the Zariski sheafification is not needed when forming the affine Grassmannian $\Gr_\cG$: in \Cref{thm:main} below, we show that the latter often agrees already with the presheaf quotient $L\cG/L^+\cG$. This appears to be new already for reductive groups $G$ over $\bC$, although for $\GL_n$ it essentially follows from \cite{Hitchin-torsors}*{Theorem 2.1.24}, and is based on the finer variant of \Cref{thm:CF-input} recorded in \Cref{thm:CF2} below, which, in addition to the geometry of $\Bun_G$, uses Quillen patching for torsors over $\bA^1_A$ to progress beyond semilocal $A$. For this variant, the relevant condition on $\cG$ is the following.

\bd[\cite{split-unramified}*{Definition 8.1}]
A semisimple group scheme $G$ over a scheme $S$ is \emph{totally isotropic} if in the canonical decomposition of \cite{SGA3IIInew}*{expos\'e XXIV, proposition 5.10~(i)}:
\[
\tst 
G^\ad \cong \prod_{i \in \{ A_n, B_n, \dotsc, G_2 \}} \Res_{S_i/S}(G_i)
\]
 of its adjoint quotient $G^\ad$, where $i$ ranges over the types of connected Dynkin diagrams, $S_i$ is a finite \'etale $S$-scheme, and $G_i$ is an adjoint semisimple $S_i$-group with simple geometric fibers of type $i$, Zariski locally on $S$ each $G_i$ has a parabolic $S_i$-subgroup that contains no $S_i$-fiber of $G_i$ (equivalently, Zariski locally on $S$ each $\Res_{S_i/S}(G_i)$ contains a nontrivial split torus $\bG_{m,\, S}$, see \cite{SGA3IIInew}*{expos\'e~XXVI, corollaire 6.12} and \cite{torsors-regular}*{end of Section 1.3.4}). 
\ed

\beg
Slightly informally, $G$ is totally isotropic if Zariski locally on $S$ it has a parabolic subgroup containing no factor of the adjoint group $G^\ad$. To see this, recall that parabolic subgroups of $G$ correspond to those of $G^\ad$, which correspond to collections of parabolic subgroups of $G_i$, one for each $i$, see \cite{torsors-regular}*{end of Section 1.3.4}.  Certainly, every quasi-split (so also every split) semisimple group is totally isotropic. 
\eeg

\bthm[\cite{totally-isotropic}*{Theorem 4.2}] \label{thm:CF2}
Let $G$ be a totally isotropic reductive group scheme over a ring $A$. For a $G$-torsor $E$ over $\bP^1_A$, if $E|_{\{t = \infty\}}$ is trivial, then $E|_{\bA^1_A}$ is also trivial.  \QED
\ethm

\bthm \label{thm:main}
For a ring $A$ and an $A\llb t \rrb$-group scheme $\cG$ that is an extension of a finite \'etale $A\llb t \rrb$-group $\cG^\et$ by a reductive $A\llb t \rrb$-group scheme $\cG^0$ that descends to a reductive $A$-group scheme whose adjoint quotient is totally isotropic, the affine Grassmannian $\Gr_\cG$ is the presheaf quotient $L\cG/L^+\cG$, that is,
\[
\Gr_\cG \cong L\cG/L^+\cG.
\]
\ethm

\bpf
As in the proof of \Cref{thm:Zariski}, by the modular description of \Cref{prop:modular} and the possibility of varying $A$, we need to show that no nontrivial $\cG$-torsor $\cE$ over $A\llb t \rrb$ trivializes over $A\llp t \rrp$. For this, as in the proof of \Cref{prop:semilocally}, \Cref{cor:L-Lplus}~\ref{m:LL-b} immediately reduces us to the case when $\cG = \cG_0$ and $\cG$ is the base change of a reductive $A$-group scheme $G$ whose adjoint quotient is totally isotropic. 

To treat this case, we again patch $\cE$ with the trivial $G$-torsor over $\bP^1_A \setminus \{ t= 0\}$ to build a $G$-torsor $E$ over $\bP^1_A$ such that $E_{A\llb t \rrb} \simeq \cE$ and $E|_{\{t = \infty\}}$ is trivial. By \Cref{thm:CF2}, this last condition forces $E|_{\bA^1_A}$ to be trivial. However, then $\cE$ is trivial, too, as desired. 
\epf

We recall from \cite{Hitchin-torsors}*{Theorems 2.1.24 and 3.1.7} that the following global variant of \Cref{cor:H1-local} was known when $G$ is either a pure inner form of $\GL_n$ or a torus. 

\bcor
For a reductive group scheme $G$ over a ring $A$ with $G^{\mathrm{ad}}$ totally isotropic, no nontrivial $G$-torsor over $A$ trivializes over $A\llp t \rrp$, in other words, we have
\[
\Ker(H^1(A, G) \ra H^1(A\llp t \rrp, G)) = \{ *\}. 
\]
\ecor

\bpf
\Cref{prop:modular} and \Cref{thm:main} show that no nontrivial $G$-torsor over $A\llb t \rrb$ trivializes over $A\llp t \rrp$. Thus, it suffices to recall from \cite{Hitchin-torsors}*{Theorem 2.1.6} that 
\[
H^1(A, G) \isomto H^1(A\llb t \rrb, G). \qedhere
\] 
\epf

Finally, we note that the preceding proof shows that the affine Grassmannian in the totally isotropic case may even be formed using polynomial loops as follows. 

\bthm \label{thm:alg-loops}
For a reductive group scheme $G$ over a ring $A$ with $G^{\mathrm{ad}}$ totally isotropic, its affine Grassmannian $\Gr_G$ may be formed as the presheaf quotient using the polynomial loops, more precisely,~we~have
\[
\Gr_G \cong L_{\mathrm{poly}} G/L^+_{\mathrm{poly}} G,
\]
explicitly, we have
\[
\ba
\tst G(A\llp t \rrp)/G(A\llb t \rrb) &\cong\tst G(A\{ t \}[\f1t])/G(A\{ t \}) \\ &\cong\tst G((A[t]_{1 + tA[t]})[\f1t])/G(A[t]_{1 + tA[t]}) \cong G(A[t, t\i])/G(A[t]),
\ea
\]
equivalently,
\[
\ba
\tst G(A\llp t \rrp) &= G(A[t, t^{-1}])G(A\llb t \rrb), \\
\tst G(A\{ t \}[\f1t]) &= G(A[t, t^{-1}])G(A\{ t \}), \\
\tst G((A[t]_{1 + tA[t]})[\f1t]) &= G(A[t, t^{-1}])G(A[t]_{1 + tA[t]}).
\ea
\]
\ethm

These equalities seem elementary, but we do not know how to argue them directly even for $G = \GL_n$. For instance, in terms of Beauville--Laszlo patching used in the proof below, one would need to argue that every finite projective $A[t]$-module that is free both over $A[t, t\i]$ and over $A\llb t \rrb$~is~free. 

\bpf
Since $G$ is affine, its functor of points preserves the Cartesianness of the squares from \Cref{cor:L-Lplus}~\ref{m:LL-a}. In particular, the map 
\[
G(A[t, t\i])/G(A[t]) \hra G(A\llp t \rrp)/G(A\llb t \rrb)
\]
is injective, and so are its counterparts for algebraic or Henselian loops in place of polynomial loops. Thus, all we need to show is that this map is also surjective or, equivalently, that
\[
G(A[t, t\i])\backslash G(A\llp t \rrp)/G(A\llb t \rrb) = \{ *\}. 
\]
However, if this double quotient was nontrivial, then we could use patching (for instance, \cite{Hitchin-torsors}*{Lemma~2.2.11~(b)}) to build a nontrivial $G$-torsor over $\bA^1_A$ that would trivialize both over $\bG_{m,\, A}$ and also over the formal completion along $\{t = 0\}$. We could then extend this $G$-torsor to all of $\bP^1_A$ by patching with the trivial torsor at infinity, and thus obtain a contradiction to \Cref{thm:CF2}. 
\epf

\brem
We do not know the extent to which the assumptions of \Cref{thm:main,thm:alg-loops} are optimal because they are imposed by our proofs, for instance, we do not know whether the affine Grassmannian $\Gr_G$ agrees with the presheaf quotient $LG/L^+G$ for every (possibly not totally isotropic) reductive $A$-group $G$, although we expect that it does not.  On the other hand, since our proofs are reductions to general results about torsors over $\bP^1_A$, the reader will have no trouble adapting them to various close variants of the affine Grassmannian that are sometimes considered in the literature, for instance, to affine Grassmannians constructed using general relative Cartier divisors in $\bA^1_A$ in place of $\{t = 0\}$. 
\erem

\subsection*{Acknowledgements} 
This article was written in relation to the Jacques Herbrand prize of the French Academy of Sciences. I thank the Academy for the prize and for the invitation to submit to the Comptes Rendus de l'Acad\'emie des Sciences on this occasion. I thank the referee for helpful remarks and suggestions. I thank Alexis Bouthier, Ofer Gabber, and Alex Youcis for helpful conversations and correspondence.  This project has received funding from the European Research Council (ERC) under the European Union's Horizon 2020 research and innovation programme (grant agreement No.~851146). This project is based upon work supported by the National Science Foundation under Grant No.~DMS-1926686.

%%%%%%%%%%%%%%%%%%%%%%%%%%%%%%%%%%

%\subfile{}

\end{document}